\documentclass[11pt,twoside]{book}

\evensidemargin \oddsidemargin
\usepackage{graphicx}
\usepackage{graphics}
\usepackage{amsfonts,amssymb, amsmath, amsthm}

\begin{document}

\begin{center}
{\large\bf{Dmitri Scheglov}}
\end{center}
\begin{center}
{\large\bf{University of Oklahoma, 2011}}
\end{center}
\begin{center}
{\large\bf{E-mail: dscheglov@math.ou.edu}}
\end{center}

\begin{center}

{\large\bf \uppercase{Canonical metrics in a conformal class}}
\end{center}

\textbf{1. Introduction and main results.} 
\

It is an old and difficult question in differential geometry if there exists a "best" or canonical metric on a smooth manifold, which makes a manifold " most symmetric". Standard examples are round spheres and flat tori, where the word " best"  means constant curvature. If there are no assumptions made about a manifold then there is a high chance that there are no reasonble " best" metrics,  which was informally explained by Gromov in [2].  Usually one assumes the existence of some geometric or algebraic structure on the manifold and considers only a class of metrics compatible with the structure.
\

Classical examples include  Kahler metrics on compact complex manifolds, left- or biinvariant metrics on Lie groups or metrics compatible with a conformal structure. And canonical metric should ideally be uniquely defined by the structure.
\

One of the first well-known  examples of canonical metrics is of course a hyperbolic metric an a compact Riemann surface of genus $g\geq 2$. In case of Kahler manifolds the Yau's proof of Calabi conjecture provides the existence of a distinguished metric in the same Kahler class as the initial one [6].
\

For a compact smooth manifold by the solution to the Yamabe problem, achieved in works of Trudinger [5], Aubin [1] and Schoen[4], each conformal structure on a compact manifold supports a metric of a constant scalar curvature. However this metric is in general not unique in the case of positive scalar curvature.
\

In a more recent work Habermann and Jost [3] construct canonical metrics in a conformal class using Green function of the Yamabe operator. Their construction requires local conformal flatness of a class if the dimension of the manifold is greater than 3.
\

In this paper we construct canonical metrics in a given conformal class for a $2n$-dimensional oriented compact smooth manifold $M$, with non-trivial $n$-th de Rham cohomology and some natural non-degeneracy assumption on the conformal class. 
\

We use Hodge theory of harmonic forms and the key point, which makes the construction very explicit is a well-known observation  that $n$-dimensional harmonic forms of a $2n$-dimensional manifold remain harmonic under conformal change of a metric.
\

We define a functional $E$ on the space of all Riemannian metrics  invariant under the natural action of the group $Diff$, which we call a Harmonic Energy.  Informally speaking the functional $E$ measures the failure of a wedge product of two harmonic forms to be harmonic. Then we prove that  inside a given conformal class there exists a unique normalized metric minimizing $E$.
\

 Moreover we obtain an explicit formula for the extremal metric in terms of the initial metric representing the conformal class and an orthonormal basis of harmonic $n$-forms. We also explicitely compute the critical Harmonic Energy and observe that corresponding value can be defined for any conformal class without any non-degeneracy assumptions and so may serve as a conformal invariant of any closed oriented smooth $2n$-manifold.
\

In the next chapter we apply our construction to Riemann surfaces, thus producing a natural family of metrics on them.
\

From now and further by manifold we mean a closed oriented smooth manifold.
\

\textbf{ Acknowledgements.} The author would like to thank Giovanni Forni for the very useful idea to consider the problem inside a fixed conformal class, without which this work would not appear.

\

\textbf{1.1. Hodge theory.} 
\

Here we briefly discuss the basics of the Hodge theory and prove a lemma which will be in use later.
\

Let  V be an oriented $d$-dimensional Euclidean space. Then the inner product produces a natural isomorphism $V=V^{*}$, which naturally extends to the graded isomorphism of Grassman algebras $\bigwedge^{k} V = \bigwedge^{k} V^{*}.$  The latter isomorphism provides a geometric interpretation of the space of $k$-forms $\Omega^{k}(V).$ Namely, for any oriented $k$-plane $P$, generated by the ordered set of orthonormal vectors $e_{1},\cdots,e_{k}$  we define a so-called decomposable $k$-form $\omega_{P}$ as follows. For a set of vectors $v_{1},\cdots,v_{k}$, the value $\omega_{P}(v_{1},\cdots,v_{k})$ equals to the algebraic $k$-volume of the projection of the parallelepiped $<v_{1},\cdots,v_{k}>$ onto $P$. In terms of the isomorphism above one can check the identity $\omega_{P}=e_{1}\wedge\cdots\wedge e_{k}$. Any other $k$-form $\omega$ is a linear combination of decomposable $k$-forms. 
 \

The inner product of two decomposable $k$-forms $\omega_{P}$, $\eta_{Q}$ is defined as a Jacobian of the orthogonal projection of oriented $k$-plane $P$ into oriented $k$-plane $Q$ and then extended by linearity to $\Omega^{k}(V)$.
\

The Hodge Star operator $*:\Omega^{k}\rightarrow\Omega^{d-k}$ is also defined in terms of decomposable forms and then extended by linearity on $\Omega^{k}(V)$. For any positively oriented orthonormal basis $e_{1},\cdots,e_{d}$ we put by definition $*(e_{1}\wedge\cdots\wedge e_{k})$=$e_{k+1}\wedge\cdots\wedge e_{d}$.
\

Now we consider a $d$-dimensional  Riemannian manifold $(M,g)$ and naturally define the Hodge Star operator $\Omega^{k}(M)\rightarrow\Omega^{d-k}(M)$ as the Riemannian metric $g$ defines inner product on the tangent space $T_{x}M$ at each point $x\in M$.
\

We then define the $L^{2}$- inner product of two $k$-forms $\omega, \eta\in\Omega^{k}(M)$ as $(\omega,\eta)_{L^{2}}=\int (\omega,\eta)_{x}Vol_{g}$ where $Vol_{g}$ is a Riemannian volume form associated to $g$ and we have a well-known equality $(\omega,\eta)_{L^{2}}=\int \omega\wedge *\eta$.
\

For a differential $d:\Omega^{k-1}\rightarrow\Omega^{k}$ there exists a formal adjoint with respect to the inner product defined above and given by the formula $\delta:\Omega^{k}\rightarrow \Omega^{k-1}$, $\delta=(-1)^{k}*^{-1}d*$.
\

The Hodge Laplacian is defined as
\

$\Delta:\Omega^{k}(M)\rightarrow\Omega^{k}(M)$
\

$\Delta=d\delta+\delta d$
\

The form $\omega$ is called harmonic if $\Delta\omega=0$ and it is easy to check that $\Delta\omega=0$ iff $d\omega=0, d*\omega=0$.
\

Let $\mathcal{H}^{k}(M)$ be the space of harmonic $k$-forms. Since all harmonic forms are closed we have a natural map $\Phi:\mathcal{H}^{k}(M)\rightarrow H^{k}(M)$ , where $H^{k}(M)$ is a de Rham cohomology of $M$ and a celebrated theorem of Hodge asserts that $\Phi$ is an isomorphism.
\

\

\textbf {Lemma 1.1.} Let $V$ be a $2n$-dimensional oriented space with an inner product $g$ and let $*_{g}:\Omega^{n}(V)\rightarrow \Omega^{n}(V)$  be a corresponding Hodge Star operator acting on $n$-forms. Let also $k$ be a positive constant. Then $*_{g}=*_{kg}$.
\

\textbf{Proof.} Obvious from the definition of the Hodge Star. \textbf{Q.E.D.}
\

\

\textbf{Corollary 1.2.} Let $M$ be a $2n$-dimensional  manifold, $g$ be a Riemannian metric, $\rho\in C^{\infty}(M), \rho>0$  and $\omega$, $\eta$ be a pair of harmonic $n$-forms with respect to $g$. Then $\omega, \eta$ are harmonic with respect to $\rho g$ and $(\omega,\eta)_{L^{2}}^{g}=(\omega,\eta)_{L^{2}}^{\rho g}$.
\

\textbf{Proof.} The first statement follows from the Lemma 1.1. and the fact that  harmonic forms $\omega$ satisfy $d\omega=0$, $d*\omega=0$. 
\

For the second statement the harmonicity of forms is not essential and it follows from the Lemma 1.1. and a formula mentioned before $(\omega,\eta)_{L^{2}}=\int\omega\wedge *\eta$. \textbf{Q.E.D.}
\

\

\textbf{2. Harmonic Energy}.
\

\

Let $M$ be a $2n$-dimensional manifold and $g$ be a Riemannian metric on it. Then all spaces $\Omega^{k}(M)$ inherit an $L^{2}$-inner product defined above and the $n$-th cohomology space $H^{n}(M)$ becomes a Euclidean space via the Hodge isomorphism $\Phi$. Then $H^{n}(M)\otimes H^{n}(M)$ naturally becomes a Euclidean space and if $e_{1}\cdots,e_{p}$ is an orthonormal basis for $H^{n}(M)$ then $e_{i}\otimes e_{j}, 1\leq i,j\leq p$ is an orthonormal basis for $H^{n}(M)\otimes H^{n}(M)$. 
\

 For cohomology classes $\left[\omega\right]$,  $\left[\eta\right]$ and $\left[\omega\wedge\eta\right]$ let $\omega_{H}$,  $\eta_{H}$, $(\omega\wedge\eta)_{H}$ be the corresponding harmonic representatives and let $\overline{\Omega^{2n}(M)}$  be a completion of $\Omega^{2n}(M)$ with respect to the $L^{2}$-norm. Then a linear operator $A:H^{n}(M)\otimes H^{n}(M)\rightarrow \overline{\Omega^{2n}(M)}$ is defined by:
\

\

 $A(\left[\omega\right]\otimes\left[\eta\right])$= $\omega_{H}\wedge\eta_{H} - (\omega\wedge\eta)_{H} $.  
\

\

 $A$  measures the failure of a wedge product of two harmonic forms to be harmonic. It is a linear operator from Euclidean to a Hilbert space and has a natural norm $\left|A\right|$, given by the formula $\left|A\right|^{2}=Tr(A^{*}A)$. Any choice of orthonormal basis $e_{1},\cdots, e_{p}\in H^{n}(M)$ allows us to write an explicit formula:
\

\

$\left|A\right|^{2}=\sum A(e_{i}\otimes e_{j})\cdot A(e_{i}\otimes e_{j})$
\

\

where $\omega\cdot\eta$ denotes the $L^{2}$-product of $2n$-forms $\omega$, $\eta$ inside $\overline{\Omega^{2n}(M)}$
\

\

\textbf{Definition 1.} Let $g$ be a Riemannian metric on $M$. Then its Harmonic Energy is
\

 $E(g)=\left|A\right|^{2}$
\

\

\textbf{Definition 2.} Let $g$ be a metric on  $M^{2n}$, such that $\int Vol_{g}=1$. Then its Normalized Conformal Class is by definition:
\

 $C(g)=\lbrace\rho g| \int\rho^{n} Vol_{g}=1, \rho\in C^{\infty}(M), \rho>0\rbrace$
\

 which is just a set of conformally equivalent metrics with total volume equal to one.
\

\

\textbf{MAIN THEOREM.} Let $(M^{2n}, g_{0})$ be a Riemannian manifold such that for any $x\in M$ there exists a form $\omega\in\mathcal{H}^{n}(M),\omega(x)\neq 0$.     
\

 Then there exists a unique metric $g\in C(g_{0})$ minimizing $E$.
 \
 
\textbf{Remark.} The assumption of the theorem of course immediately implies $H^{n}(M)\neq 0$.
\

\

 \textbf{PROOF.}
 \
 
 Let us pick a set of harmonic forms $\xi_{1},\cdots,\xi_{p}$ which form an orthonormal basis in $H^{n}(M)$ for any metric from $C(g_{0})$. This is a crucial ingredient of the proof and such a choice is possible because of the Corollary 1.2.
 \
 
 Now we consider any metric $g=\rho g_{0}\in C(g_{0})$  and calculate $E(g)$ using that $g$-harmonic $2n$-forms are proportional to the $Vol_{g}$:
\

 $E(g)=\sum A(\xi_{i}\otimes\xi_{j})\cdot A(\xi_{i}\otimes\xi_{j})$ $=$ $\sum\int A(\xi_{i}\otimes\xi_{j})\wedge*A(\xi_{i}\otimes\xi_{j})$ $=$ $\sum\int(\xi_{i}\wedge\xi_{j}-(\int\xi_{i}\wedge\xi_{j})Vol_{g})\wedge*(\xi_{i}\wedge\xi_{j}-(\int\xi_{i}\wedge\xi_{j})Vol_{g})$
 \
 
 Let us now introduce the following notations: let $f_{ij}$ be smooth functions, defined by $\xi_{i}\wedge\xi_{j}=f_{ij}Vol_{0}$ and $c_{ij}$ be constants, defined by $c_{ij}=\int\xi_{i}\wedge\xi_{j}$.
 \
 
 Using that $Vol_{g}=\rho^{n} g_{0}$, we can rewrite
\

 $E(g)=\sum\int (f_{ij}\rho^{-n}-c_{ij})^{2}\rho^{n}Vol_{0}$. Opening brackets we obtain 
\

$E(g)=\int f^{2}\rho^{-n}Vol_{0} -C^{2}$ where $f=\sqrt{\sum f_{ij}^{2}}$ and $C=\sqrt{\sum c_{ij}^{2}}$. 
 \
 
 Now let us prove by contradiction that $f(x)>0$ for any $x\in M$. Indeed  $f(x)=0$$\Longrightarrow$$ f_{ij}(x)=0$$\Longrightarrow $ $\omega\wedge\eta (x)=0$ for any  $\omega,\eta\in\mathcal{H}^{n}(M)$ as $i,j$ run through the basis of $\mathcal{H}^{n}(M)$. But then for any $\omega\in\mathcal{H}^{n}(M)$, $(\omega,\omega)_{x}Vol$=$\omega\wedge *\omega(x)=0$ which implies $\omega(x)=0$ and contradicts the theorem assumption.
\

To find a minimum of $E(g)$ among all metrics $g\in C(g_{0})$ we use a form of integral Cauchi inequality on the functions $f, \rho$:
\

$\int\frac{f^{2}}{\rho^{n}}Vol_{0}\cdot\int\rho^{n}Vol_{0}\geq\left(\int f Vol_{0}\right)^{2}$
\

which immediately implies that the extremal metric is $g=\rho g_{0}$ with
\

 $\rho=(f/\int f Vol_{0})^{1/n}$  and minimal $E(g)=(\int f Vol_{0})^{2}-C^{2}$. \textbf{Q.E.D.}
\

\

\textbf{Remark 1.} Using explicit formulas from the proof one can easily check that the minimizing metric and critical harmonic energy are indeed independent on the initial metric $g_{0}$ we started with.
\

\

\textbf{Remark 2.} As the expression for the critical energy $E(g)$ is independent on the initial metric $g_{0}$ one can define the harmonic energy of a given conformal class on any compact smooth oriented manifold $M^{2n}$ dropping the assumption of the theorem, which is only required to guarantee the non-degeneracy of the critical metric.
\

\

\textbf {3. Canonical metrics on Riemann surfaces.}
\

Consider a closed oriented surface $M$ of genus $g\geq 2$. The conformal classes of metrics on $M$ are in natural one-to-one correspondence with complex structures on it. Let us fix a conformal class on $M$. Then harmonic 1-forms for such a conformal class are precisely the real parts of abelian differentials for the corresponding complex structure. 
\

As it is well known that for any $x\in M$ there exists an abelian differential $w, w(x)\neq 0$ we have that either $Re(w)\neq 0$ or $Im(w)=Re(-i w)\neq 0$ which means that natural conformal class of any Riemann surface satisfies the assumptions of the main theorem.
\ 

So we produced a canonical metric for any Riemann surface.
\

\

\textbf{4. Open questions.}
\

\textbf{1.} It seems to be an interesting question to describe the properties of the metric on the Riemann surface, defined above, and in particular to see if its curvature is negative on $M$.
\

\textbf{2.} A manifold is called formal if there exists a metric such that the wedge product of any two harmonic forms is harmonic. It is not hard to show that any closed surface of genus $g\geq 2$ is not formal which implies that the value of critical Harmonic Energy is a positive smooth function $E_{g}$ on the moduli space $\mathcal{M}_{g}$ of complex curves of genus {g}. 
\

As $\mathcal{M}_{g}$ is not compact it is an interesting question if $E_{g}$ is strictly positive on $\mathcal{M}_{g}$ which can be reformulated as if $M$ is formal "at infinity" of the moduli space.

\end{document}